# Moment-based approximation of cell and face integrals in THINC/QQ reconstruction on unstructured grids

Young-Lin Yoo[1,*]

[1]SEMO Simulation, Paju, Republic of Korea

* Correspondence: younglin90@gmail.com

**Abstract**

On unstructured grids, the tangent of hyperbola for interface capturing (THINC) method with quadratic surface representation and Gaussian quadrature (THINC/QQ) requires topology-dependent numerical quadrature for cell and face integrals and Newton iteration for the surface constant. This study introduces a quadratic moment-matched sigmoid (QMMS) approximation that evaluates these quantities from the mean and variance of the quadratic field while retaining the hyperbolic-tangent reconstruction profile and the boundary variation diminishing (BVD) formulation. A slope-matched Gaussian cumulative distribution function and a moment-matched Gaussian model provide closed-form approximations for the surface constant and face averages, eliminating Newton iteration and topology-dependent runtime quadrature from these evaluations. Single-cell tests and nine two- and three-dimensional benchmarks assess reconstruction accuracy, flow-field differences, and computational cost. Across the tested cases, the principal transported and shock structures remain closely aligned with those obtained using THINC/QQ, with local differences concentrated mainly near discontinuities, contact regions, and developed shear layers. Reconstruction-stage profiling gives mean THINC/QQ-to-QMMS time ratios of 2.12 in two dimensions and 1.85 in three dimensions. QMMS therefore reduces reconstruction-stage cost while replacing topology-dependent runtime integral evaluation with a common moment-based formulation.

## 1. Introduction

The tangent of hyperbola for interface capturing (THINC) method was introduced as an algebraic volume-of-fluid (VOF) reconstruction based on a bounded hyperbolic-tangent profile, with the profile location determined from the prescribed cell average [1]. Early multidimensional developments included THINC/WLIC [2], orientation-dependent interface scaling [3], and continuous multidimensional reconstruction for curved interfaces [4]. THINC was subsequently incorporated as a discontinuity-capturing candidate in the boundary variation diminishing (BVD) scheme, which selects among reconstruction candidates according to the variation across cell boundaries [5]. Higher-order BVD formulations combined polynomial reconstructions with THINC candidates of different steepness [6,7], and related hybrid formulations incorporated targeted essentially non-oscillatory (TENO) and weighted essentially non-oscillatory (WENO) reconstructions [8–10]. In compressible multiphase flow, THINC was also developed as an algebraic interface-sharpening procedure [11] and incorporated into BVD reconstruction to resolve smooth flow structures and moving material interfaces with low numerical dissipation [12].

On unstructured grids, multidimensional THINC reconstruction requires cell- and face-dependent integration. Continuous THINC reconstruction was extended to triangular and tetrahedral cells [13], followed by quadratic-surface representations in two and three dimensions [14]. THINC with quadratic surface representation and Gaussian quadrature (THINC/QQ) combines a multidimensional quadratic function with Gaussian quadrature on hybrid unstructured grids [15]. Subsequent refinements addressed coordinate normalization, polynomial-coefficient evaluation, higher-order representations, quadrature, and time discretization [16]. In this formulation, determining the surface constant requires numerical integration of the nonlinear reconstruction profile, and evaluating face averages requires additional quadrature. For mixed-element grids, the baseline implementation uses topology-dependent mappings and integration rules. THINC/QQ has been combined with BVD reconstruction for compressible flow on unstructured grids [17]. Related THINC-based formulations have been applied to interfacial flows, three-dimensional unstructured grids, adaptive and moving grids, overset grids, phase-changing flows, and multi-material hydrodynamics [18–24].

Related developments address interface representation and reconstruction cost. Coupling THINC with level-set information enabled conservative high-order surface representations [25], while THINC scaling extended polynomial interface representations to structured and unstructured grids [26]. Quartic surface reconstruction further increased the geometric representation order on unstructured grids [27]. Direct temporal integration has been used to reduce the time-integration cost of THINC [28], while neural-network-assisted BVD selection reduces the number of candidate reconstructions [29]. THINC/HAM replaces Newton–Raphson iteration for the interface position with a homotopy-analysis-based iterative algorithm [30]. Alternative sigmoid functions provide another approach to algebraic interface reconstruction. The sigmoid functions for interface capturing (SFINC) framework extends THINC to other smooth sigmoid functions and considers exact or approximate one-dimensional interface-location relations [31]. Related VOF methods use piecewise logistic profiles on unstructured grids [32] or steepness-adjustable harmonic functions [33]. The admissible steepness of THINC-type reconstruction has also been studied [34], and THINC/WLIC developments have examined geometric fidelity and interface sharpness separately [35]. Outside the algebraic THINC family, piecewise-paraboloid reconstruction uses an iterative geometric procedure on arbitrary three-dimensional grids [36].

This study proposes a quadratic moment-matched sigmoid (QMMS) approximation to the cell and face integrals in THINC/QQ while retaining the hyperbolic-tangent profile. A slope-matched Gaussian cumulative distribution function (CDF) surrogate and a Gaussian model of the within-cell distribution of the quadratic field yield an approximate cell-average relation based on the field mean and variance. Algebraic inversion of this relation provides the surface constant, and the same closure supplies approximate face averages. The resulting expressions use precomputed geometric moments and require no Newton iteration or runtime cell and face quadrature for reconstruction. QMMS is compared with THINC/QQ within the same BVD scheme, with the monotonic upstream-centered scheme for conservation laws (MUSCL) providing the polynomial reconstruction candidate. Single-cell tests cover triangular, quadrilateral, tetrahedral, hexahedral, prismatic, and pyramidal cells, while flow benchmarks on triangular and tetrahedral grids evaluate accuracy, discontinuity resolution, and reconstruction cost.

## 2. Numerical methods

### 2.1. Baseline scheme: MUSCL–THINC/QQ–BVD

The scalar advection equation and the compressible Euler equations are discretized using a cell-centered finite-volume method on unstructured grids. The advected scalar is reconstructed for scalar advection, and the primitive variables are reconstructed for the Euler equations. The numerical flux is specified for each benchmark, and time integration uses a third-order strong-stability-preserving (SSP) Runge–Kutta method [37]. Three candidates are considered for each reconstructed variable. The first is a MUSCL reconstruction with the Barth–Jespersen limiter [38] and local bounds obtained from cells sharing a vertex with the target cell [39,40]. The other two are THINC/QQ reconstructions with steepness parameters $\beta_s = 0.8$ and $\beta_l = 1.4$ [17]. The BVD principle [5,17] selects among these candidates using the total boundary variation (TBV), evaluated as the sum of the absolute jumps in the reconstructed states across the interior faces of the cell.

The THINC/QQ reconstruction in cell $\Omega_i$ [15,17] is expressed as

$$q_i(\boldsymbol{x}) = \bar{q}_i^{min} + \frac{\bar{q}_i^{max} - \bar{q}_i^{min}}{2}\left[1 + \sigma\left(\hat{\beta}_i P_i(\boldsymbol{x}) + D_i\right)\right], \tag{1}$$

where $\bar{q}_i^{min}$ and $\bar{q}_i^{max}$ are the minimum and maximum cell-average values over the vertex-neighbor stencil of cell $\Omega_i$. The quadratic polynomial $P_i(\boldsymbol{x})$ defines the shape of the local reconstruction surface. The scaled steepness is $\hat{\beta}_i = \beta/h_i$, where $h_i$ is the characteristic cell length. The sigmoid is denoted by $\sigma$, with $\sigma = \tanh$ for THINC/QQ. With $\bar{q}_i$ denoting the prescribed cell average, the dimensionless surface constant $D_i$ is determined from the cell-average constraint

$$\frac{1}{|\Omega_i|}\int_{\Omega_i} \sigma\left(\hat{\beta}_i P_i(\boldsymbol{x}) + D_i\right) d\Omega = \frac{2\left(\bar{q}_i - \bar{q}_i^{min}\right)}{\bar{q}_i^{max} - \bar{q}_i^{min}} - 1, \tag{2}$$

where $|\Omega_i|$ denotes the cell area in two dimensions and the cell volume in three dimensions. In the baseline scheme, the surface constant $D_i$ is determined by Gaussian quadrature and Newton iteration [15].

The baseline implementation uses the tanh addition identity to reduce repeated transcendental-function evaluations [15]. The terms $\tanh\left(\hat{\beta}_i P_i\right)$ are precomputed at the quadrature points, and Newton iteration is then performed on the transformed unknown $\tanh(D_i)$ using algebraic operations. Each residual evaluation requires a weighted sum over the quadrature points. Once $D_i$ has been determined, the face average on the jth face $\Gamma_{ij}$ of cell $\Omega_i$, denoted by $\bar{q}_{ij}^f$, is defined as

$$\bar{q}_{ij}^f = \frac{1}{\left|\Gamma_{ij}\right|}\int_{\Gamma_{ij}} q_i\,(\boldsymbol{x})\, d\Gamma, \tag{3}$$

where $\left|\Gamma_{ij}\right|$ denotes the edge length in two dimensions and the face area in three dimensions.

The face average $\bar{q}_{ij}^f$ is used to evaluate the numerical flux. The baseline implementation uses topology-dependent Gaussian quadrature for the cell integral in Eq. (2) and the face average in Eq. (3). The cell rules use 6 and 9 points for triangular and quadrilateral cells in two dimensions and 11, 12, 12, and 16 points for tetrahedral, pyramidal, prismatic, and hexahedral cells in three dimensions, respectively [15]. Face quadrature uses 4 points per edge in two dimensions and 6 and 9 points on triangular and quadrilateral faces in three dimensions, respectively. Newton iteration starts from $D_i = 0$ and terminates when the absolute residual of the quadrature-discretized cell-average constraint is less than $10^{-8}$.

### 2.2. Moment-based approximation of the cell and face integrals

QMMS retains the quadratic representation, candidate steepness parameters, and BVD selection procedure of the baseline scheme described in Section 2.1 and modifies only the evaluation of the surface constant and face averages. Specifically, the numerical quadrature and Newton iteration used in the baseline formulation are replaced by closed-form moment-based approximations. To construct these approximations, QMMS uses a slope-matched Gaussian-CDF surrogate $\sigma_G$, defined by

$$\sigma_G(s) = \mathrm{erf}\left(\frac{\sqrt{\pi}}{2}s\right) = 2\Phi\left(\sqrt{\frac{\pi}{2}}s\right) - 1, \tag{4}$$

where erf is the error function and $\Phi$ is the standard normal CDF. The scaling of $\sigma_G(s)$ gives $\sigma_G'(0) = 1$, matching the slope of $\tanh(s)$ at the origin. For the quadratic field $P_i(\boldsymbol{x})$, the cell mean $m_i$ and variance $v_i$ are

$$m_i = \frac{1}{|\Omega_i|}\int_{\Omega_i} P_i(\boldsymbol{x})\, d\Omega,\ v_i = \frac{1}{|\Omega_i|}\int_{\Omega_i} [P_i(\boldsymbol{x}) - m_i]^2\, d\Omega. \tag{5}$$

For a general quadratic field $P_i(\boldsymbol{x})$, evaluating $m_i$ and $v_i$ requires geometric moments through second and fourth order, respectively. The geometric moments are precomputed analytically in two dimensions and by Gaussian quadrature after a Duffy transformation in three dimensions [41], whereas the field mean and variance are updated with the reconstruction coefficients.

QMMS approximates the generally non-Gaussian distribution of $P_i(\boldsymbol{x})$ under uniform spatial weighting within the cell by a Gaussian with mean $m_i$ and variance $v_i$. Related approximations for sigmoid expectations under Gaussian inputs are discussed in [42,43]. For a Gaussian random variable $s$ with mean $m_i$ and variance $v_i$, the convolution identity [44] gives

$$\mathbb{E}[\sigma_G(s)] = \sigma_G\left(\frac{m_i}{\sqrt{1 + \frac{\pi}{2}v_i}}\right), \tag{6}$$

where $\mathbb{E}[\cdot]$ denotes expectation with respect to $s$. Applying Eq. (6) to the affine argument $\hat{\beta}_i s + D_i$ and replacing the resulting surrogate sigmoid by tanh gives

$$\frac{1}{|\Omega_i|}\int_{\Omega_i} \tanh\left[\hat{\beta}_i P_i(\boldsymbol{x}) + D_i\right] d\Omega \approx \tanh\left[\frac{\hat{\beta}_i m_i + D_i}{\sqrt{1 + \frac{\pi}{2}\hat{\beta}_i^2 v_i}}\right]. \tag{7}$$

The resulting relation remains approximate even for Gaussian inputs because the surrogate sigmoid is replaced by tanh. Section 3.1 assesses the combined effects of the distribution and sigmoid approximations. With $\bar{Q}_i$ denoting the right-hand side of Eq. (2), algebraic inversion gives

$$D_i = \sqrt{1 + \frac{\pi}{2}\hat{\beta}_i^2 v_i}\, \mathrm{atanh}(\bar{Q}_i) - \hat{\beta}_i m_i. \tag{8}$$

For normalized cell averages strictly between −1 and 1, this closed-form relation gives the approximate surface constant without Newton iteration or runtime cell quadrature. The face average in Eq. (3) is approximated by applying the closure in Eq. (7) over $\Gamma_{ij}$, with $m_i$ and $v_i$ replaced by the face mean $m_{ij}^f$ and variance $v_{ij}^f$, respectively, and $D_i$ retained from Eq. (8). These moments are defined as in Eq. (5), with $\Omega_i$ replaced by $\Gamma_{ij}$, and are evaluated from precomputed face geometric moments.

## 3. Numerical results

Single-cell tests and nine flow benchmarks are used to compare QMMS with THINC/QQ. For each flow benchmark, the two methods use the same grid, Courant–Friedrichs–Lewy (CFL) number, numerical flux, boundary conditions, and time-integration scheme. Unless otherwise stated, spatial coordinates, time, and flow variables are reported in nondimensional form.

### 3.1. Single-cell reconstruction accuracy

The tests considered triangular and quadrilateral cells in two dimensions and tetrahedral, hexahedral, prismatic, and pyramidal cells in three dimensions, including distorted configurations. In two dimensions, each vertex was generated by uniformly sampling its angular direction over a full revolution about the cell center and its radial distance over 0.25–1. A generated cell was retained only if it was convex, the angular separation between adjacent vertices exceeded 0.45 rad, every edge length exceeded 0.12, and the cell area exceeded 0.06. In three dimensions, the vertices of reference tetrahedral, hexahedral, prismatic, and pyramidal cells were randomly perturbed by 5–45% of the reference edge length. A generated cell was retained only if it remained non-folded and star-shaped with respect to its centroid and its volume exceeded 15% of the reference volume.

For the quadratic field, the linear direction was sampled isotropically. In two dimensions, the orientation angle was drawn uniformly over a full revolution, and the three quadratic coefficients were sampled independently from a uniform distribution on a symmetric interval, with a common half-width drawn uniformly between 0 and 2 and divided by the characteristic cell length. In three dimensions, the linear direction was generated as an isotropic random unit vector, and the six Hessian inputs were sampled independently from a standard normal distribution, multiplied by a common factor drawn uniformly between 0 and 2, and divided by the characteristic cell length before conversion to the quadratic coefficients. The prescribed normalized cell average was sampled with equal probability from either a uniform distribution over $[-0.999, 0.999]$ or the hyperbolic tangent of a variable sampled uniformly over $[-4, 4]$. The surface-constant and face-average evaluation stages were examined separately, with $10^6$ samples generated for each cell shape. The steepness parameter was sampled uniformly over $0.5 \le \beta \le 5.0$, a range encompassing $\beta_s = 0.8$ and $\beta_l = 1.4$ used in the benchmark calculations.

Reference values used to assess reconstruction accuracy were obtained by high-order numerical reintegration of the tanh profile with the surface constant obtained by each method. To apply the reference quadrature to cells of different shapes, the cells and their faces were decomposed into simplices. The reference quadrature used $128 \times 128$ points per triangular subcell and 128 points per edge in two dimensions, and $16 \times 16 \times 16$ points per tetrahedral subcell and $16 \times 16$ points per triangular subface in three dimensions. For surface-constant evaluation, the cell-average error $\varepsilon_D$ is the absolute difference between the reintegrated normalized cell average and the prescribed value in Eq. (2). The face-average error $\varepsilon_\mathrm{F}$ is the absolute difference between the computed face average and the corresponding high-order reference face average, evaluated using the same surface constant.

To analyze the cell-average error associated with surface-constant evaluation in the single-cell tests, the scaled-field magnitude $s_{\max}$ is defined as

$$s_{\max} = \hat{\beta}_i \max_{x\in\Omega_i} |P_i(\boldsymbol{x})|. \tag{9}$$

This quantity combines the quadratic-field amplitude and the steepness parameter into a single measure. The value of $s_{\max}$ measures the maximum magnitude of the scaled quadratic field, allowing samples with different field amplitudes and steepness parameters to be compared on a common scale.

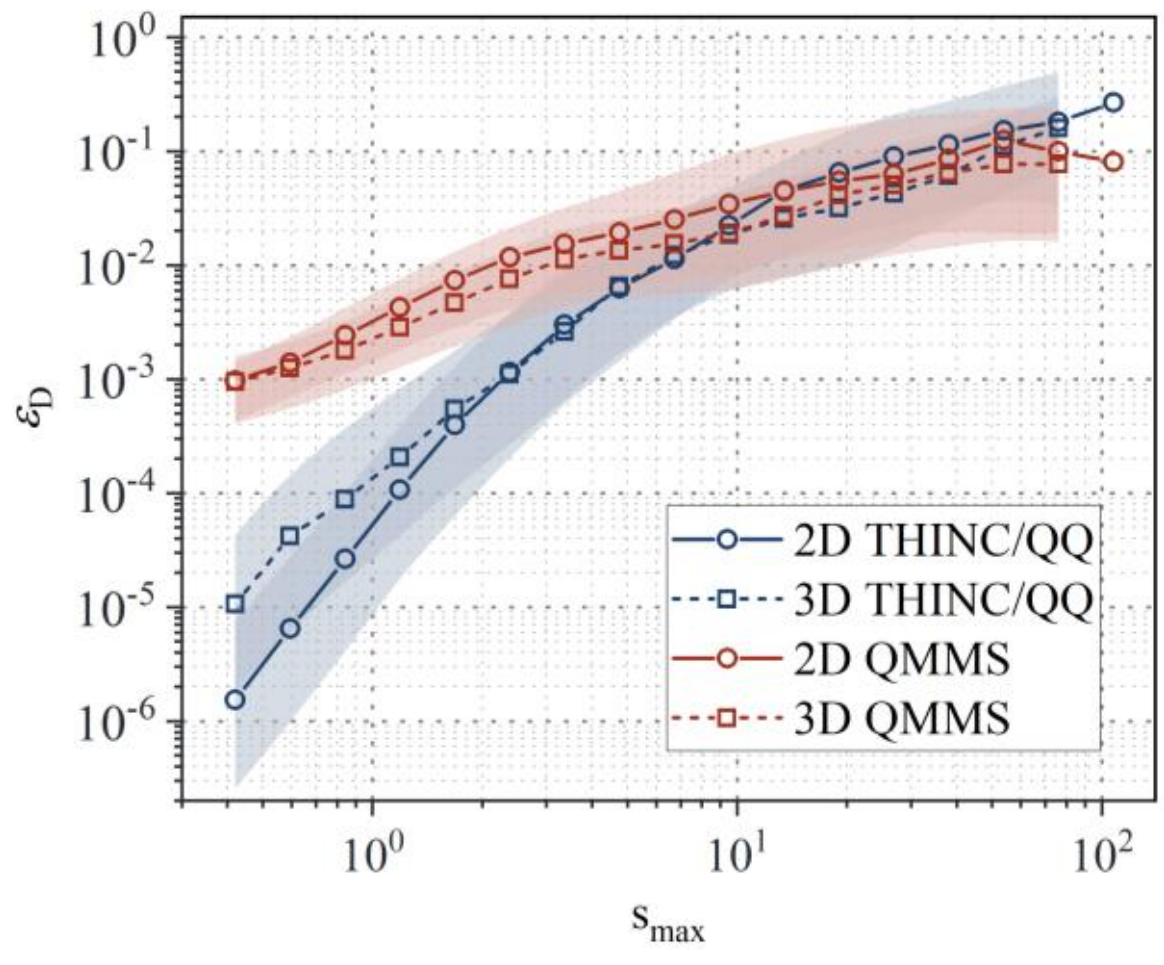

Fig. 1. Median cell-average error $\varepsilon_D$ with respect to the scaled-field magnitude $s_{\max}$.

Figure 1 shows the median cell-average error $\varepsilon_D$ in each $s_{\max}$ bin. The symbols and connecting lines denote the median values, while the shaded bands span the 25th–75th percentiles within each bin. For $s_{\max} < 10$, THINC/QQ generally gives lower median errors than QMMS. At larger $s_{\max}$, the median errors of the two reconstructions approach the same order of magnitude. The increase in the THINC/QQ error with $s_{\max}$ becomes less steep toward the upper end of the logarithmic plot. For THINC/QQ, the surface constant satisfies the quadrature-discretized constraint to the Newton tolerance, whereas the reported error is evaluated by high-order reintegration of the resulting profile. QMMS instead evaluates the cell-average approximation from the moments $m_i$ and $v_i$ rather than by pointwise quadrature, and its median error increases more gradually. At larger $s_{\max}$, the QMMS interquartile range also spans a wider absolute error interval. Samples with similar $s_{\max}$ can have different spatial distributions of $P_i(x)$, and distributional features beyond the mean and variance of $P_i$ are not retained by the moment approximation. These unresolved differences can contribute to the spread in the errors. The two- and three-dimensional results exhibit the same overall trends.

Table 1 summarizes the reconstruction accuracy and computational costs for $\beta \leq 1.5$, including the steepness parameters used in the flow benchmarks. Table 1(a) reports the unweighted root-mean-square (RMS) errors over the sampled cells. The RMS ranges of $\varepsilon_D$ and $\varepsilon_{\mathrm{F}}$ overlap between THINC/QQ and QMMS, but the upper limits of the face-average RMS ranges are larger for QMMS in both dimensions. Table 1(b) reports the computational time for the corresponding single-cell evaluations. For surface-constant evaluation, the algebraic moment-based closure in QMMS provides speedups of 21.0–21.7× in two dimensions and 19.8–42.0× in three dimensions relative to THINC/QQ. For face-average evaluation, the corresponding speedups are 1.69× in two dimensions and 1.71–4.99× in three dimensions. The reduction in computational cost is most pronounced for surface-constant evaluation, particularly for three-dimensional cells, where the THINC/QQ formulation requires more quadrature operations.

Table 1. Cell-average errors ($\varepsilon_{\mathrm{D}}$), face-average errors ($\varepsilon_{\mathrm{F}}$), and single-cell computational costs for $\beta \leq 1.5$.

(a) RMS

| Dimension | Quantity | THINC/QQ RMS | QMMS RMS |
|---|---|---|---|
| 2D | $\varepsilon_{\mathrm{D}}$ | $2.7\times10^{-3}$–$3.1\times10^{-2}$ | $1.3\times10^{-2}$–$4.3\times10^{-2}$ |
| | $\varepsilon_{\mathrm{F}}$ | $1.4\times10^{-3}$–$1.3\times10^{-2}$ | $1.1\times10^{-2}$–$4.5\times10^{-2}$ |
| 3D | $\varepsilon_{\mathrm{D}}$ | $3.7\times10^{-3}$–$3.3\times10^{-2}$ | $4.1\times10^{-3}$–$3.0\times10^{-2}$ |
| | $\varepsilon_{\mathrm{F}}$ | $3.1\times10^{-3}$–$2.1\times10^{-2}$ | $5.6\times10^{-3}$–$4.1\times10^{-2}$ |

(b) Computational cost

| Dimension | Stage | THINC/QQ cost | QMMS cost | Speedup |
|---|---|---|---|---|
| 2D | Surface constant | 317–323 ns | 14.9–15.1 ns | 21.0–21.7× |
| | Face average | 111 ns | 65.7 ns | 1.69× |
| 3D | Surface constant | 0.524–1.16 μs | 0.0265–0.0276 μs | 19.8–42.0× |
| | Face average | 0.137–0.400 μs | 0.0765–0.0813 μs | 1.71–4.99× |

### 3.2. Two-dimensional benchmark problems

Six two-dimensional benchmarks are considered. Rigid-body rotation solves the scalar advection equation, and the other five cases solve the compressible Euler equations for a calorically perfect gas with $\gamma = 1.4$. The baseline uses the reconstruction procedure in Section 2.1 with steepness parameters $\beta_s$ and $\beta_l$.

### 3.2.1. Rigid-body rotation

The rigid-body rotation problem [45] tests the transport of continuous and discontinuous scalar profiles over one revolution. The scalar advection equation

$$\frac{\partial \phi}{\partial t} + \nabla \cdot (a\phi) = 0 \tag{10}$$

is solved in the computational domain $[0,1]^2$ with the divergence-free velocity field

$$a(x,y) = \left(-2\pi\left(y - \frac{1}{2}\right), 2\pi\left(x - \frac{1}{2}\right)\right), \tag{11}$$

which gives a rotation period of unity. The solution is advanced to $t = 1$, when the exact solution returns to the initial distribution. The initial scalar field $\phi_0$ contains a slotted cylinder centered at $(0.5, 0.75)$, a cone centered at $(0.5, 0.25)$, and a cosine hump centered at $(0.25, 0.5)$, all with radius $r_0 = 0.15$. For each object, $r$ is the distance from its center. The slotted cylinder has $\phi_0 = 1$ for $r \le r_0$, except in the slot defined by $|x - 0.5| < 0.025$ and $y < 0.85$, where $\phi_0 = 0$. The cone and cosine hump have $\phi_0 = 1 - r/r_0$ and $\phi_0 = \frac{1}{4}[1 + \cos(\pi r/r_0)]$, respectively, for $r \le r_0$, and $\phi_0 = 0$ elsewhere. A criss-cross triangulation with 200 subdivisions in each coordinate direction gives 160,000 triangular cells. An upwind flux is used with a CFL number of 0.3.

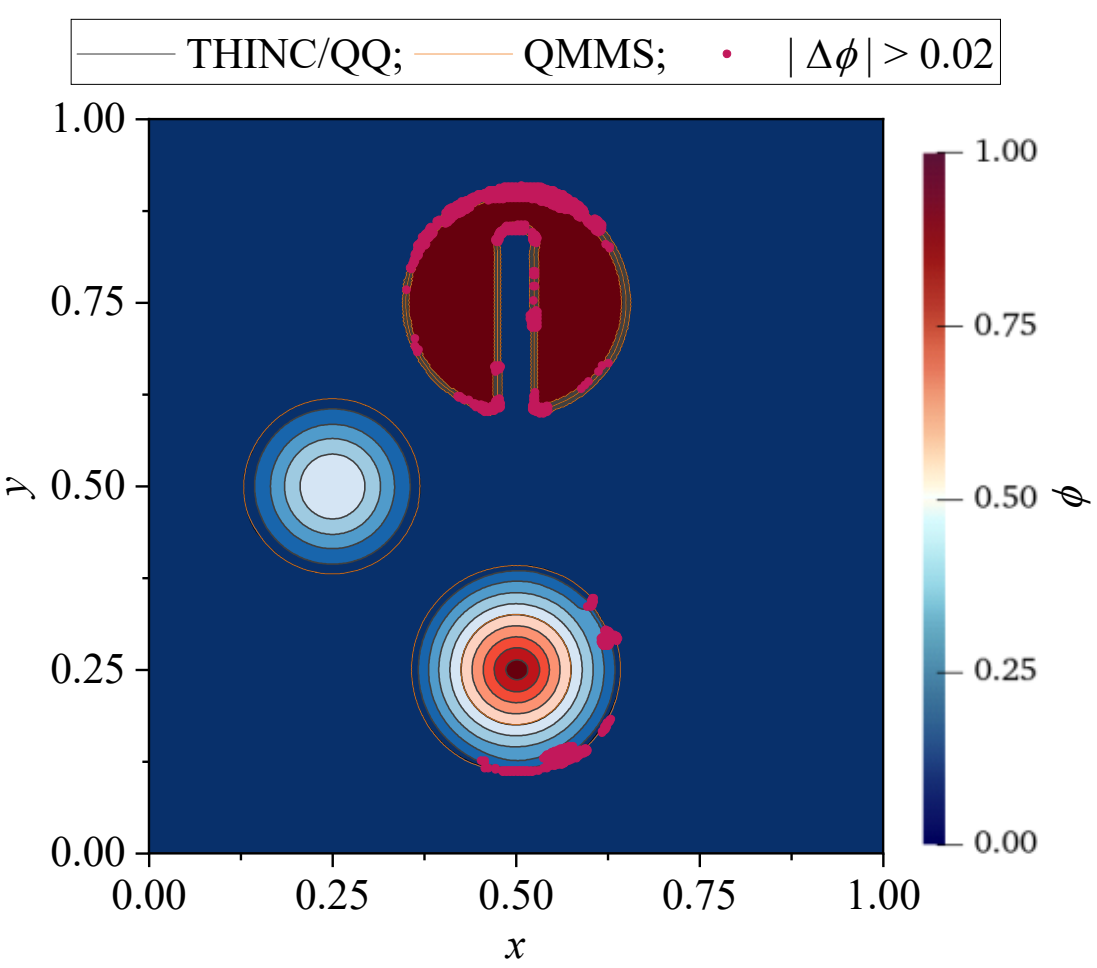

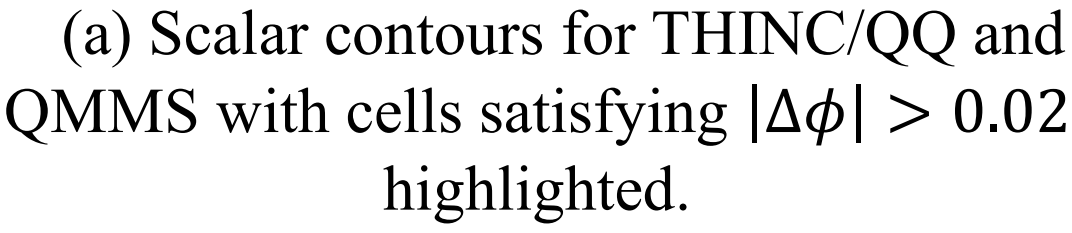

(a) Scalar contours for THINC/QQ and QMMS with cells satisfying $|\Delta\phi| > 0.02$ highlighted.

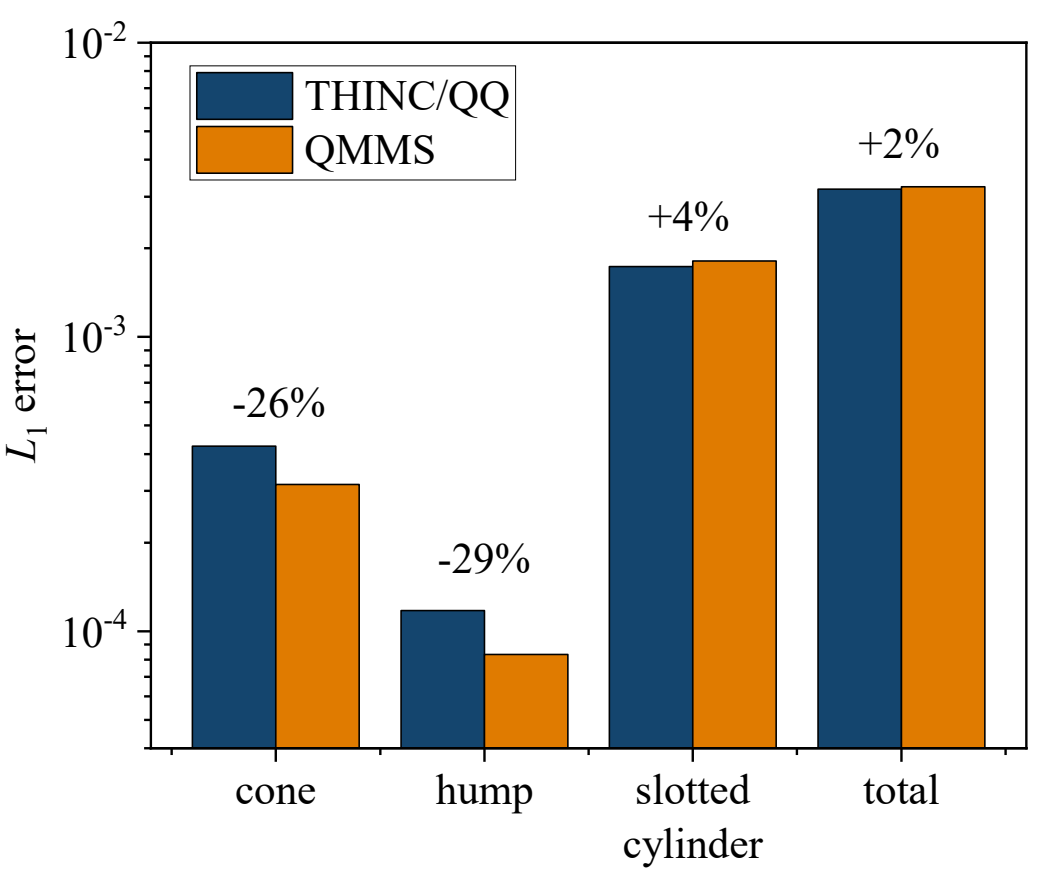

(b) Global and object-wise $L_1$ errors.

Fig. 2. Comparison of rigid-body rotation results obtained with THINC/QQ and QMMS after one revolution.

The scalar fields after one revolution show nearly overlapping contours for all three objects in Fig. 2(a). Magenta marks cells where the absolute scalar difference between QMMS and THINC/QQ exceeds 0.02. These cells account for 0.68% of the grid and occur mainly near the slotted-cylinder and cone boundaries, including the edges of the slot. The $L_1$ errors shown in Fig. 2(b) are computed as cell-area-weighted sums of absolute differences from the initialized scalar field. For the object-wise comparison, the same quantity is evaluated over cells whose centroids lie within each object's disk, including the slot for the slotted cylinder, without normalization by the disk area. The global $L_1$ errors are $3.18 \times 10^{-3}$ for THINC/QQ and $3.24 \times 10^{-3}$ for QMMS, corresponding to an increase of approximately 2%. Relative to THINC/QQ, the cone and cosine-hump errors decrease by 26% and 29%, respectively, whereas the slotted-cylinder error increases by approximately 4%.

### 3.2.2. Configuration 3 Riemann problem

The four-quadrant Configuration 3 Riemann problem [46,47] is solved on $[0,1]^2$ to examine interacting shocks and slip lines. The initial states use exact Rankine–Hugoniot relations rather than rounded tabulated values. The lines $x = 0.8$ and $y = 0.8$ partition the domain into four quadrants with constant states $(\rho, u, v, p)$:

$$\begin{gathered} Q_1 : (3/2, 0, 0, 3/2), \\ Q_2 : (33/62, u_S, 0, 3/10), \\ Q_3 : (77/558, u_S, u_S, 9/310), \\ Q_4 : (33/62, 0, u_S, 3/10), \end{gathered} \tag{12}$$

where $u_S = 4/\sqrt{11}$. Transmissive boundary conditions are imposed on all four sides. A rising-diagonal triangulation with $N = 280$ subdivisions in each coordinate direction gives 156,800 triangular cells. The local Lax–Friedrichs flux [48] is used with a CFL number of 0.4.

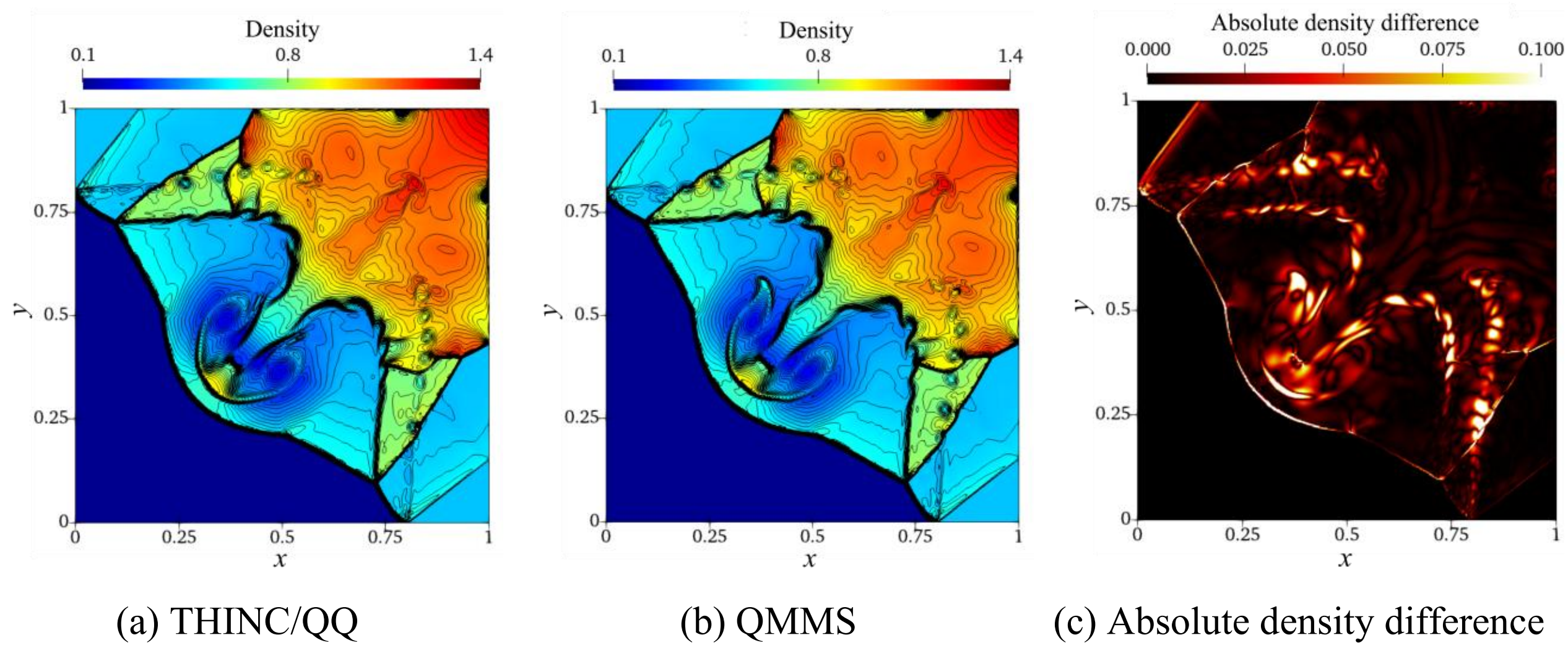

(a) THINC/QQ (b) QMMS (c) Absolute density difference

Fig. 3. Comparison of density contours obtained with THINC/QQ and QMMS for the Configuration 3 Riemann problem at $t = 0.8$.

Figures 3(a) and 3(b) compare the density fields at $t = 0.8$, and Fig. 3(c) shows the absolute density difference. Both reconstructions retain the primary shock pattern and the slip lines that roll into the central region. The local differences in Fig. 3(c) are concentrated along the shocks and the rolled shear layers. A pronounced difference appears in the central roll, where the density contours in Figs. 3(a) and 3(b) also differ in symmetry. QMMS evaluates the surface constant and face averages using algebraic moment-based approximations, which can modify the face values and the corresponding TBV values. When competing reconstructions have similar TBV values, these changes can alter BVD candidate selection and provide one possible numerical source of the observed differences. The larger differences along the rolled shear layers are consistent with the cumulative effect of these local reconstruction differences as the slip lines evolve, while the principal shock structures remain closely aligned.

### 3.2.3. Shock–vortex interaction

The shock–vortex interaction problem [48] consists of an isentropic vortex upstream of a stationary Mach 1.1 normal shock. The computational domain is $[0, 2] \times [0, 1]$, and the shock is initially at $x = 0.5$. The vortex has strength $\varepsilon = 0.3$, core radius $r_c = 0.05$, and center $(x_v, y_v) = (0.25, 0.5)$. Defining

$$r^2 = \frac{(x - x_v)^2 + (y - y_v)^2}{r_c^2},\ \tau = \varepsilon \exp\left[\frac{1}{2}(1 - r^2)\right], \tag{13}$$

the velocity and thermodynamic fields are

$$\begin{gathered} u = 1.1\sqrt{\gamma} + \tau\frac{y - y_v}{r_c},\ v = -\tau\frac{x - x_v}{r_c}, \\ \delta T = -\frac{(\gamma - 1)\varepsilon^2}{2\gamma}\exp(1 - r^2), \\ \rho = (1 + \delta T)^{1/(\gamma-1)},\ p = (1 + \delta T)^{\gamma/(\gamma-1)}. \end{gathered} \tag{14}$$

The uniform upstream state is $(\rho, u, v, p) = (1, 1.1\sqrt{\gamma}, 0, 1)$, and the corresponding downstream state from the Rankine–Hugoniot relations is approximately $(1.169, 1.114, 0, 1.245)$. The upstream state is imposed at the left boundary, and the other boundaries are transmissive. A uniform rising-diagonal triangulation with $400 \times 200$ subdivisions gives 160,000 triangular cells. The rotated Harten–Lax–van Leer–contact (HLLC) flux [49] is used with a CFL number of 0.3.

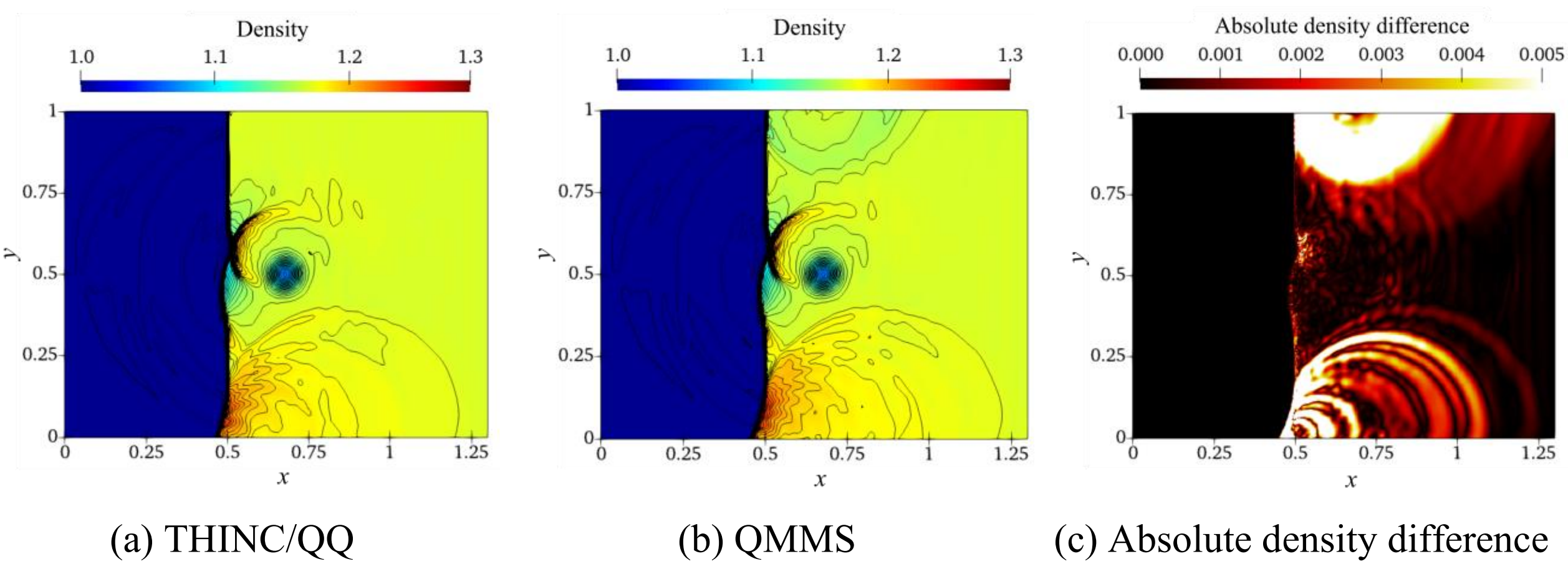

(a) THINC/QQ (b) QMMS (c) Absolute density difference

Fig. 4. Comparison of density contours obtained with THINC/QQ and QMMS for the shock–vortex interaction at $t = 0.35$.

Figures 4(a) and 4(b) show similar shock locations and vortex deformation at $t = 0.35$. In both fields, the vortex core lies just downstream of the shock, and the shock bends locally around the interaction region. Curved wave fronts extend away from the vortex toward the upper and lower parts of the domain. The density fluctuations associated with these weak acoustic waves are of order $10^{-3}$ relative to the background state. The absolute density difference in Fig. 4(c) is concentrated mainly along these waves rather than at the dominant shock or vortex core. This spatial pattern is consistent with small differences in wave position or amplitude while the dominant shock and vortex structures remain nearly aligned.

### 3.2.4. Shock–mixing-layer interaction

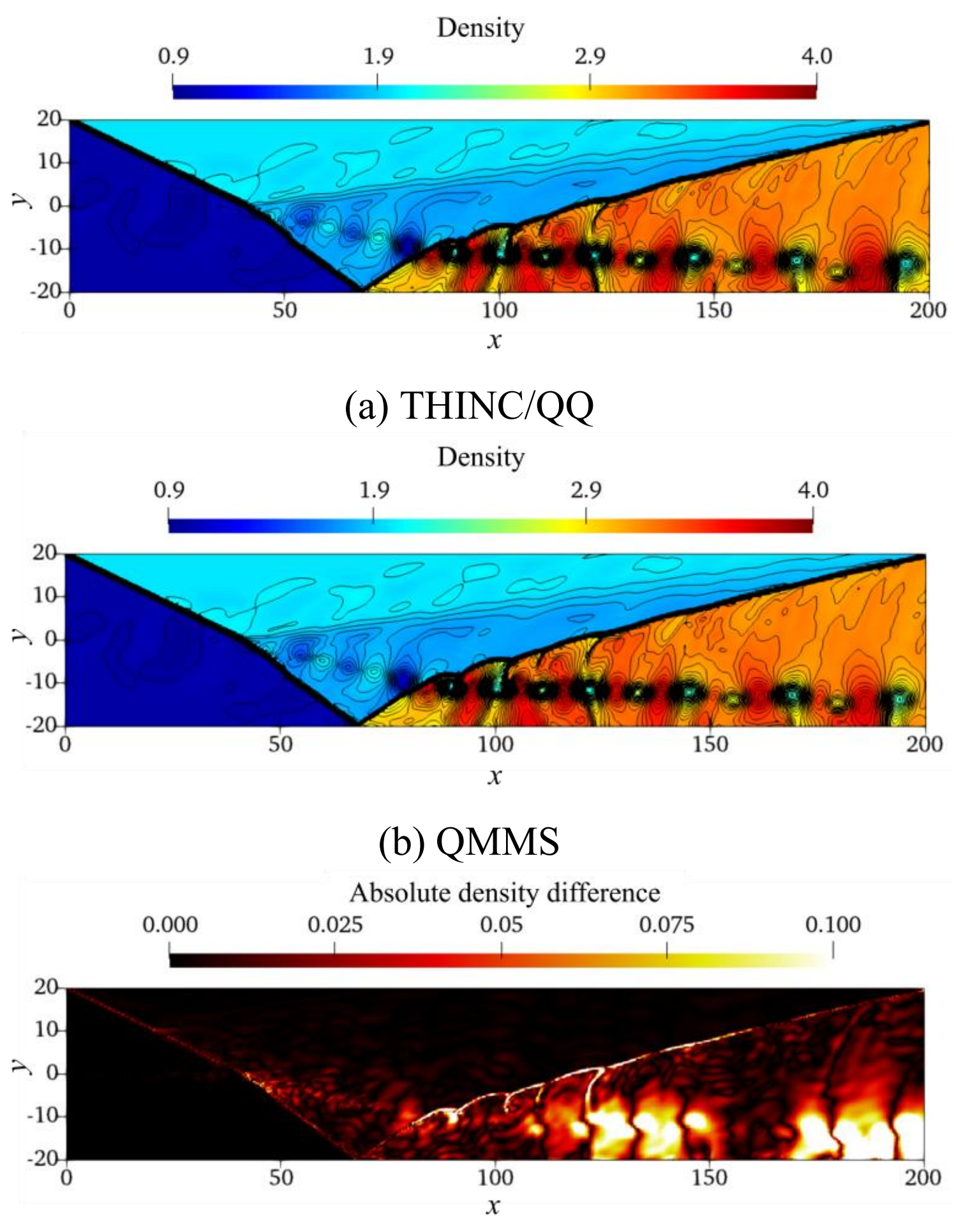

(a) THINC/QQ

(b) QMMS

(c) Absolute density difference

Fig. 5. Comparison of density contours obtained with THINC/QQ and QMMS for the shock–mixing-layer interaction at $t = 120$.

The shock–mixing-layer interaction problem is adapted from Yee et al. [50] and solved on $[0, 200] \times [-20, 20]$. An oblique shock interacts with a spatially developing shear layer and the resulting Kelvin–Helmholtz vortices. The initial state is

$$\rho = 1,\ u(y) = 2.5 + 0.5\tanh(2y),\ v = 0,\ p = 0.496. \tag{15}$$

The corresponding free-stream velocities above and below the shear layer are $u_1 = 3$ and $u_2 = 2$, respectively, with a convective Mach number of $M_c = 0.6$. Unlike the viscous reference problem, the present computation is inviscid. At the inflow, a transverse velocity perturbation is superimposed on the mean velocity profile,

$$v'(y,t) = \exp\left(-\frac{y^2}{b}\right)\sum_{k=1}^{2} a_k \cos\left(\frac{2\pi k t}{T} + \phi_k\right), \tag{16}$$

where $a_1 = a_2 = 0.05$, $\phi_1 = 0$, $\phi_2 = \pi/2$, $b = 10$, and $T = 11.194$. The post-shock state corresponding to $M_1 = 3.6$ and a 12° flow deflection is determined from the weak branch of the oblique-shock relations and prescribed at the upper boundary. A slip-wall condition is imposed at the lower boundary, while the right boundary is treated as transmissive. The computational domain is discretized using a uniform rising-diagonal triangulation with $400 \times 80$ subdivisions, yielding 64,000 triangular cells with $h_x = h_y = 0.5$. The rotated HLLC flux is used with a CFL number of 0.4.

At $t = 120$, Figs. 5(a) and 5(b) show similar density fields, including the incident oblique shock, its reflection, and the developing vortex street. The principal fronts and the sequence of rolled structures occupy similar locations at the displayed scale. The absolute density difference in Fig. 5(c) remains small upstream and becomes more pronounced within the downstream mixing layer and its wake. These differences are concentrated around the rolled structures. QMMS modifies the reconstructed face values used in the numerical flux and TBV evaluation, **so** small local differences may persist and develop as the shear layer is advected downstream.

### 3.2.5. Mach 3 forward-facing step

The Mach 3 forward-facing step of Woodward and Colella [51] tests interacting shocks and a downstream slip line. The computational domain is the channel $[0, 3] \times [0, 1]$ with a step occupying $[0.6, 3] \times [0, 0.2]$. The initial flow is the uniform Mach 3 state

$$(\rho, u, v, p) = (1.4, 3.0, 0, 1.0). \tag{17}$$

The inflow state is imposed at $x = 0$, while the right boundary is treated as transmissive. Slip-wall boundary conditions are applied to the upper and lower channel walls and the step surfaces. The computational domain is discretized using $480 \times 160$ subdivisions of the full rectangular domain with alternating diagonal splits, followed by removal of the cells within the step, yielding 129,024 triangular cells with $h = 1/160$. The Harten–Lax–van Leer (HLL) flux [52] is used with a CFL number of 0.3.

The bow shock, its upper-wall reflection, the Mach stem, and the triple point occupy comparable locations in the THINC/QQ and QMMS results at $t = 4$, as shown in Figs. 6(a) and 6(b). Differences become more visible downstream along the rolled slip line. Figure 6(c) likewise shows that the larger density differences occur within and around these rolled structures, whereas differences near the strong shocks occupy narrow bands. Several successive rolls contribute to the downstream patches. As in the shock–mixing-layer interaction, the concentration of the differences along the rolled shear layer is consistent with small reconstruction-dependent changes persisting as these structures evolve downstream.

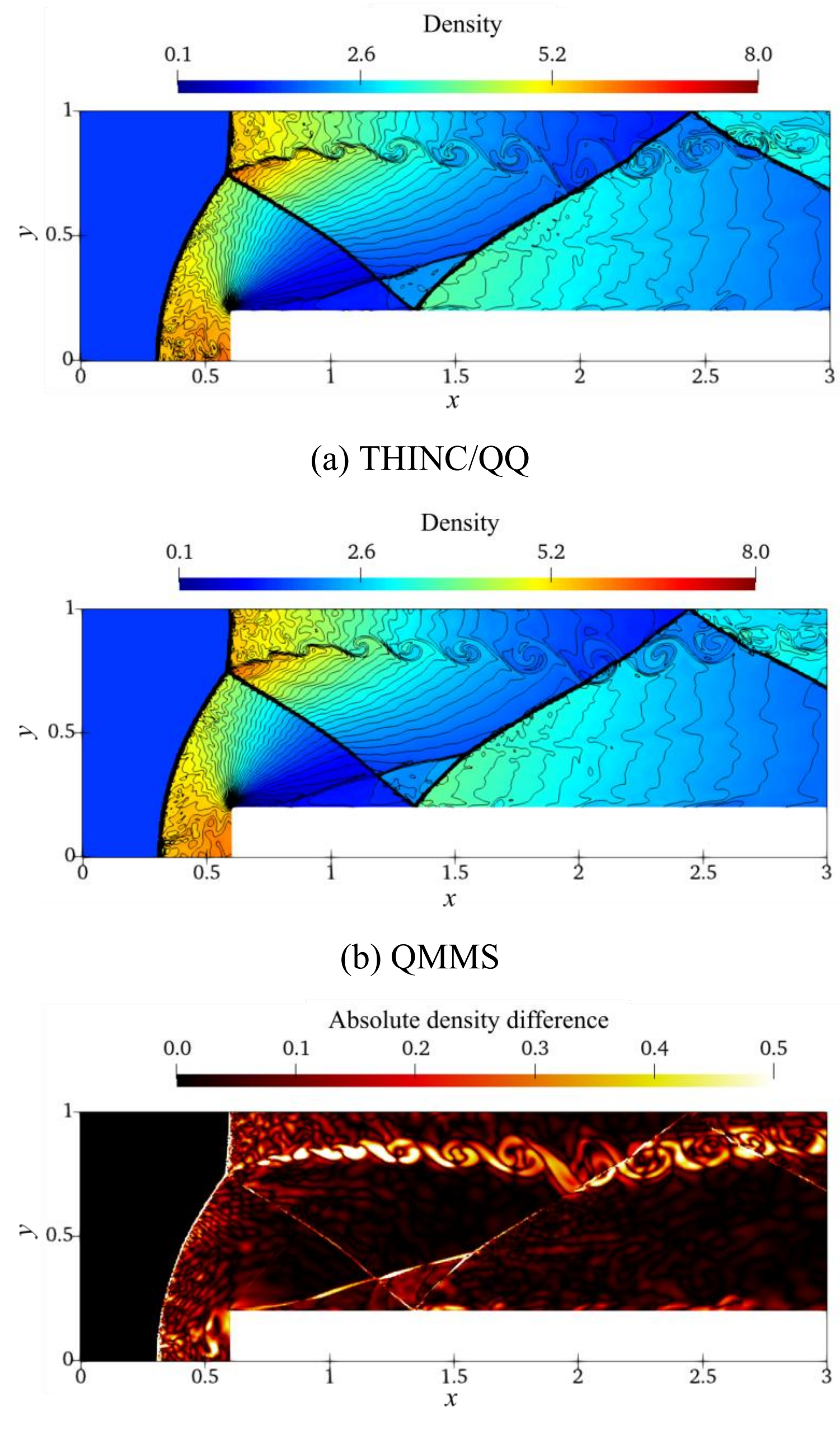

(c) Absolute density difference

Fig. 6. Comparison of density contours obtained with THINC/QQ and QMMS for the Mach 3 forward-facing step at $t = 4$.

### 3.2.6. Double Mach reflection

The Mach 10 double Mach reflection problem of Woodward and Colella [51] is solved on $[0, 4] \times [0, 1]$. Initially, the shock is inclined at 60° and intersects the lower boundary at $x = 1/6$. The pre-shock and post-shock states are

$$\begin{aligned} (\rho, u, v, p)_{\text{pre}} &= (1.4, 0, 0, 1), \\ (\rho, u, v, p)_{\text{post}} &= (8, 7.1447096, -4.125, 116.5), \end{aligned} \tag{18}$$

respectively, with the post-shock state assigned for

$$x < \frac{1}{6} + \frac{y}{\sqrt{3}}. \tag{19}$$

The post-shock state is prescribed at the left boundary and along the lower boundary for $x < 1/6$, while a slip-wall condition is imposed on the lower boundary for $x \geq 1/6$. The right boundary is treated as transmissive. At the upper boundary, the instantaneous position of the incident shock is specified by

$$x_s(t) = \frac{1}{6} + \frac{1 + 20t}{\sqrt{3}}, \tag{20}$$

with the corresponding pre-shock and post-shock states prescribed on either side of the shock. The computational domain is discretized using an approximately uniform unstructured triangulation consisting of 370,421 triangular cells and 186,212 nodes, with a nominal edge length of approximately 1/200. The rotated HLLC flux is used with a CFL number of 0.3.

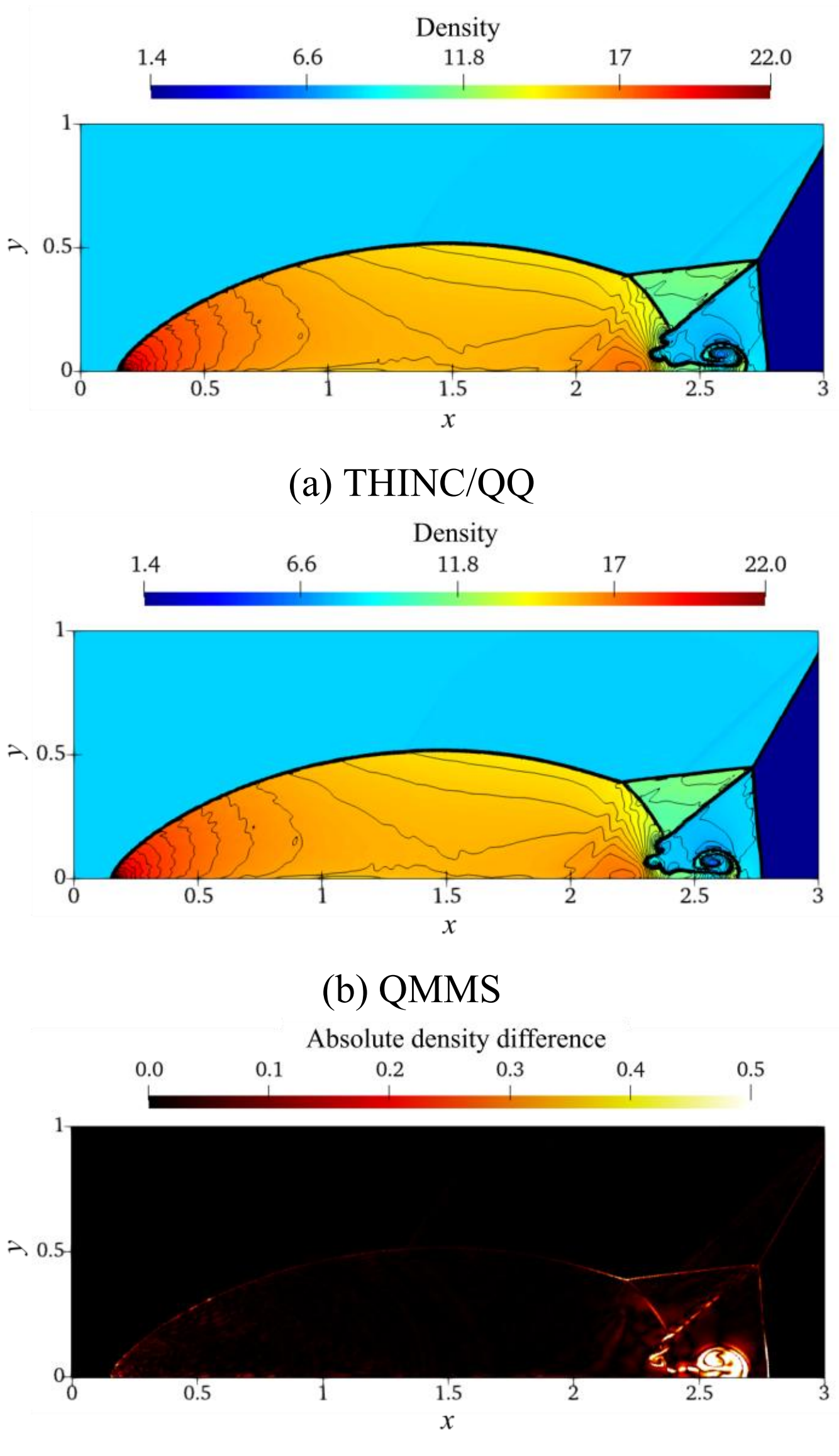

(a) THINC/QQ

(b) QMMS

(c) Absolute density difference

Fig. 7. Comparison of density contours obtained with THINC/QQ and QMMS for the double Mach reflection problem at $t = 0.2$.

The incident and reflected shock systems obtained with THINC/QQ and QMMS are closely aligned at $t = 0.2$, with similar Mach-stem and triple-point configurations, as shown in Figs. 7(a) and 7(b). The curved reflected front extends across much of the displayed domain and meets the incident-shock system above the near-wall jet. Near the lower boundary, both calculations also exhibit a compact rolled structure connected to the downstream slip lines, although local differences are visible in its internal contours. The absolute density difference in Fig. 7(c) is most pronounced around this rolled structure and the associated downstream slip lines, while narrower differences trace portions of the reflected shock front. In contrast, much of the compressed region behind the incident shock exhibits smaller differences at the displayed scale. Overall, the discrepancies remain localized primarily around the strongly deformed roll and slip-line structures, consistent with the preceding shear-layer benchmarks.

### 3.3. Three-dimensional benchmark problems

The three-dimensional benchmarks are time-reversed deformation, the spherical Riemann problem, and shock–cylinder interaction. The first solves the scalar advection equation, while the other two solve the compressible Euler equations for a calorically perfect gas with $\gamma = 1.4$. All three use tetrahedral grids and the reconstruction and time-integration procedures described in Section 2.1.

#### 3.3.1. Time-reversed deformation

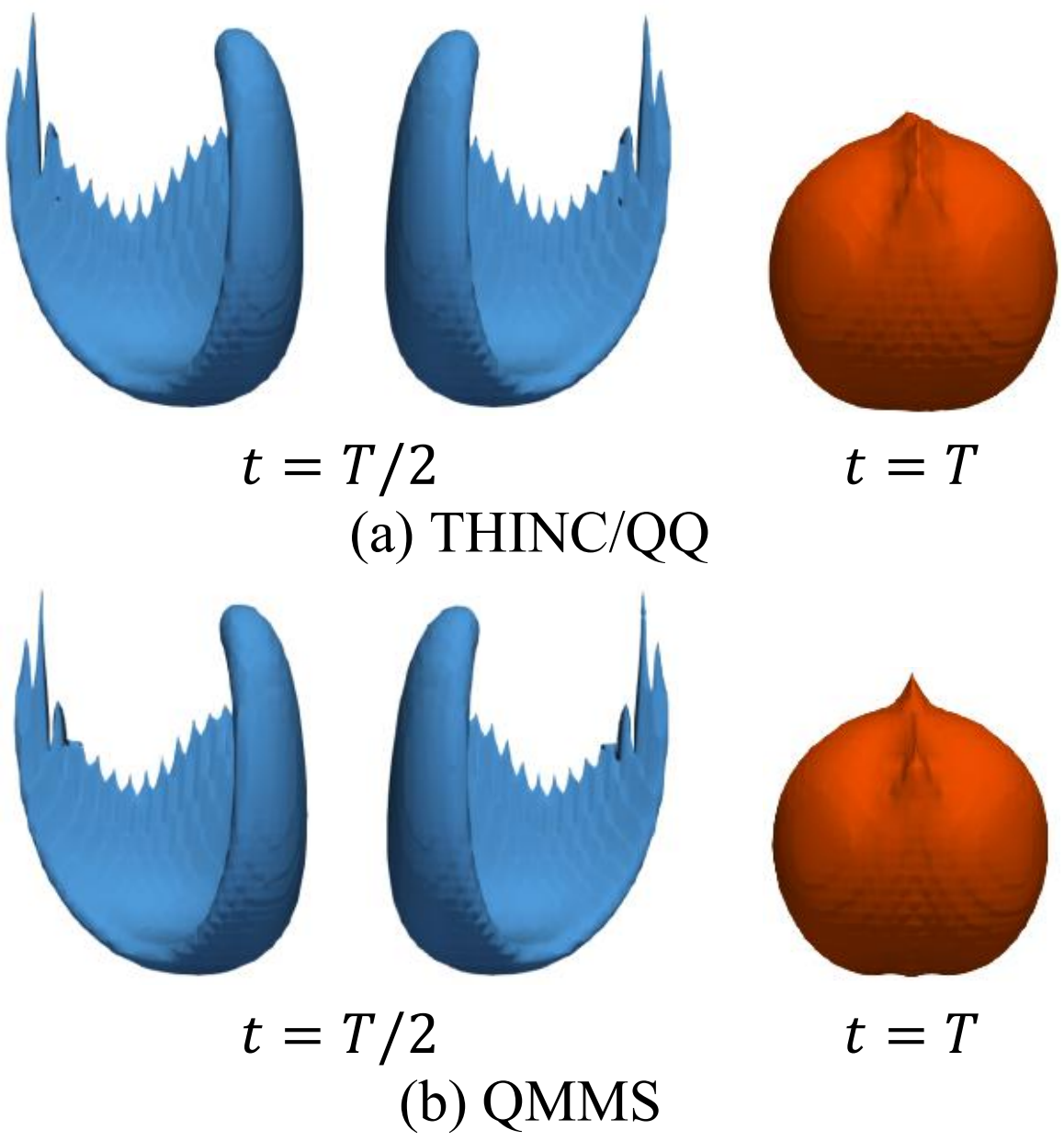

Fig. 8. Comparison of the $\phi = 0.5$ isosurfaces obtained with THINC/QQ and QMMS at $t = T/2$ and $t = T$.

The time-reversed deformation problem of Enright et al. [53] tests interface transport under stretching and folding. A sphere of radius $0.15$, initially centered at $(0.35, 0.35, 0.35)$, is advected in the unit cube by the divergence-free velocity field

$$a_x = 2\sin^2(\pi x)\sin(2\pi y)\sin(2\pi z)\cos(\pi t/T),$$
$$a_y = -\sin^2(\pi y)\sin(2\pi x)\sin(2\pi z)\cos(\pi t/T), \quad (21)$$
$$a_z = -\sin^2(\pi z)\sin(2\pi x)\sin(2\pi y)\cos(\pi t/T),$$

with $T = 3$. The velocity reverses at $t = T/2$, and the exact solution returns to the initial configuration at $t = T$. The initial scalar is $\phi = 1$ inside the sphere and $\phi = 0$ elsewhere. The grid contains 2,633,856 tetrahedral cells generated from a $76^3$ hexahedral background grid. An upwind flux is used with a CFL number of 0.25.

Figure 8 compares the $\phi = 0.5$ isosurfaces at $t = T/2$ and $t = T$. At the maximum-deformation time $t = T/2$, both reconstructions produce stretched and folded sheet-like interfaces with comparable thick regions and thin, irregular outer portions. At $t = T$, both interfaces return to an approximately spherical shape with an upper protrusion and residual surface corrugation. The $L_1$ errors relative to the initialized scalar field are $4.39 \times 10^{-3}$ for THINC/QQ and $4.58 \times 10^{-3}$ for QMMS, respectively, corresponding to an approximately 4% higher error for QMMS.

### 3.3.2. Spherical Riemann problem

The spherical Riemann problem of Langseth and LeVeque [54], also used in related numerical studies [55–57], tests shock propagation and reflection in three dimensions. The computational domain is $[0, 1.5] \times [0, 1.5] \times [0, 1]$. The gas is initially at rest with $(\rho, p) = (1, 1)$ outside a sphere of radius 0.2 centered at $(0, 0, 0.4)$, and $(\rho, p) = (1, 5)$ inside it. Symmetry conditions are imposed at $x = 0$ and $y = 0$, reflecting-wall conditions at $z = 0$ and $z = 1$, and transmissive conditions at $x = 1.5$ and $y = 1.5$. The tetrahedral grid contains 4,020,000 cells. The HLLC flux [58] is used with a CFL number of 0.4.

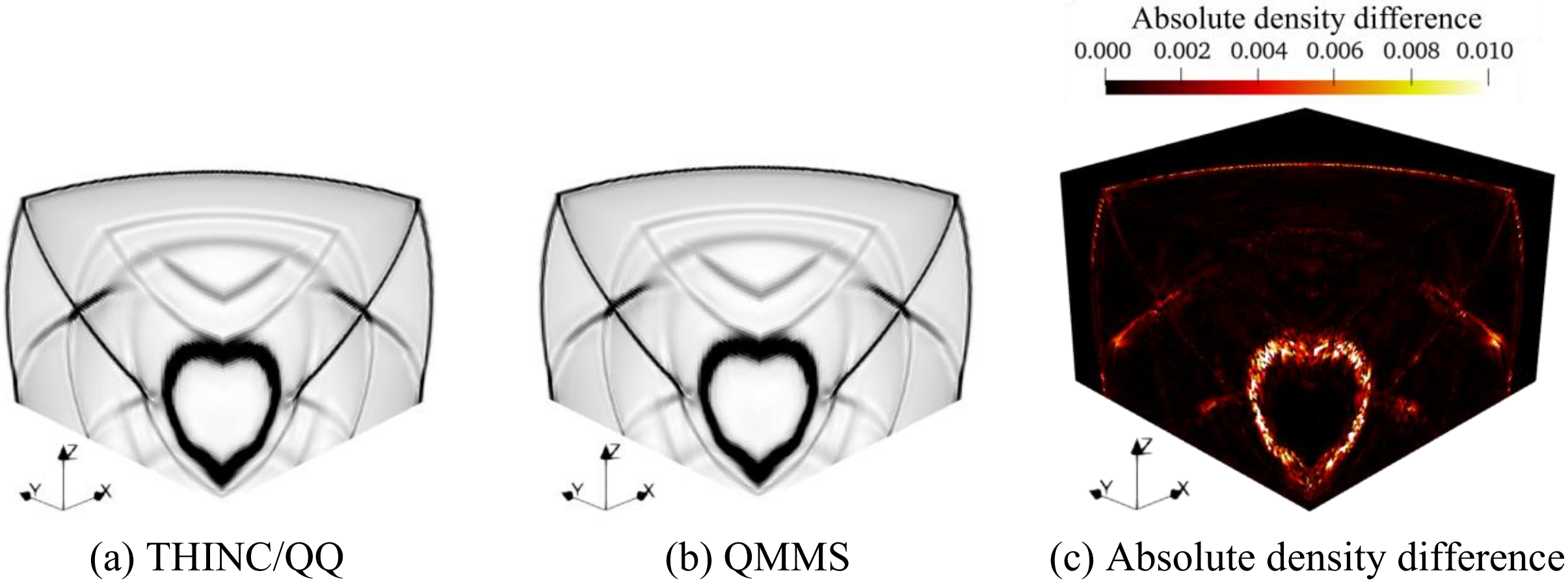

(a) THINC/QQ (b) QMMS (c) Absolute density difference

Fig. 9. Numerical schlieren fields and absolute density difference for the spherical Riemann problem at $t = 0.8$.

The numerical schlieren fields in Figs. 9(a) and 9(b) place the outward-propagating shock, reflected waves, and Mach stem at similar locations at $t = 0.8$. The density differences in Fig. 9(c) follow thin bands along these fronts and form a stronger patch around the inner curved wave structure, while much of the remaining slice shows smaller differences at the displayed scale. Minimum pressures are 0.8288 and 0.8287, and density ranges are [0.3024, 1.2533] and [0.3023, 1.2540] for THINC/QQ and QMMS, respectively. The minimum pressures differ by $10^{-4}$, and the corresponding density extrema differ by no more than $7 \times 10^{-4}$.

### 3.3.3. Shock–cylinder interaction

The shock-induced vorticity problem of Langseth and LeVeque [54], also studied by Huang et al. [59], is used to assess shock–density-interface interaction and baroclinic vorticity generation. The computational domain is $[0, 1.5] \times [0, 1] \times [0, 0.5]$, and the ambient gas is initially at rest with $(\rho, p) = (1, 1)$. Two orthogonal cylinders of radius 0.2 are embedded in the domain. The first is aligned with the $z$-axis at $x = y = 0$ and has $(\rho, p) = (1, 10)$. The second is aligned with the $y$-axis at $x = 0.4$ and $z = 0$ and has $(\rho, p) = (0.1, 1)$. Symmetry conditions are imposed at $x = 0$, $y = 0$, and $z = 0$, and the remaining boundaries are transmissive. The tetrahedral grid contains 1,864,800 cells. The rotated HLLC flux is used with a CFL number of 0.4.

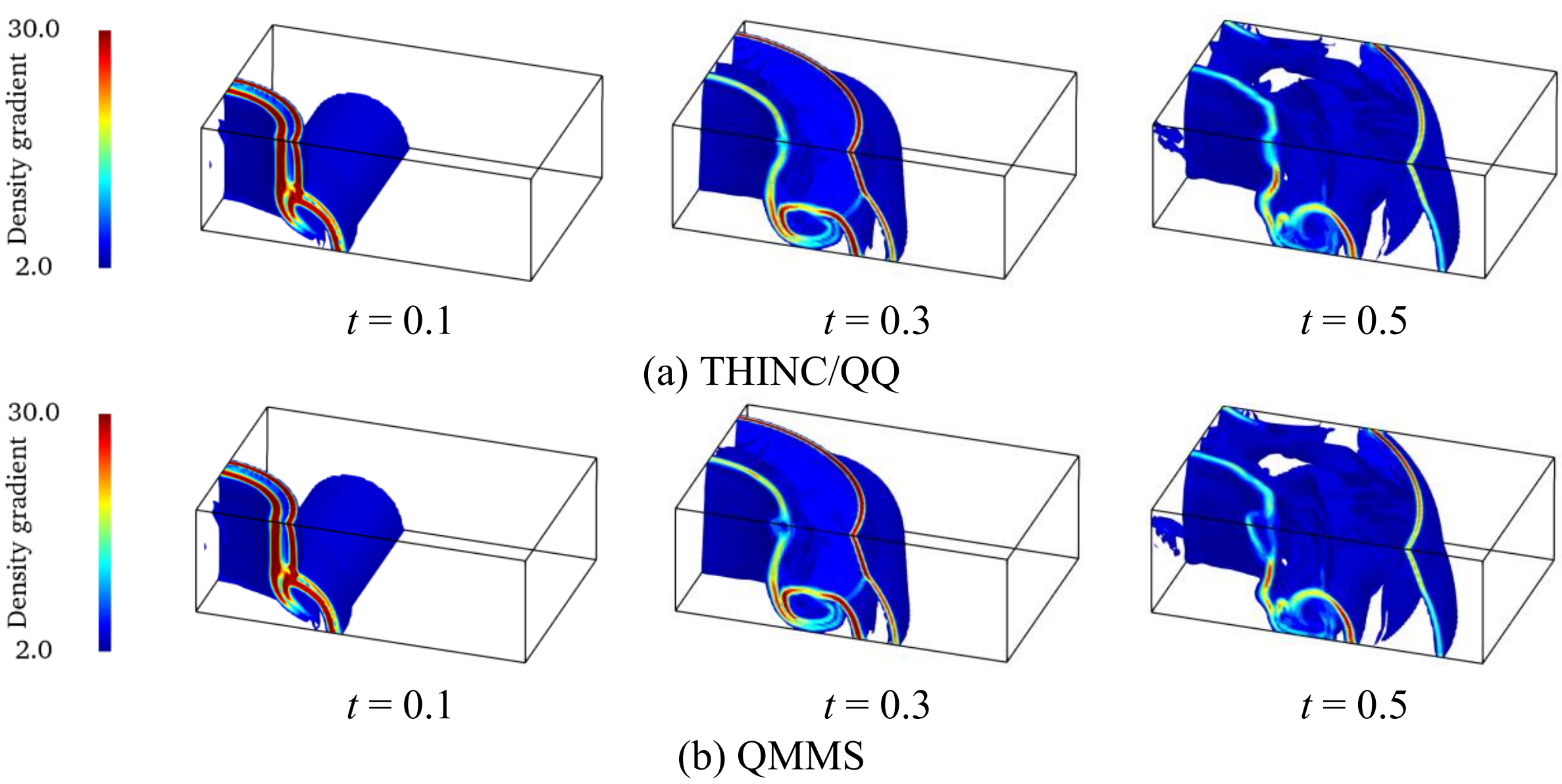

Fig. 10. Density-gradient magnitude and isosurfaces for shock–cylinder interaction at $t = 0.1, 0.3, 0.5$.

The density-gradient fields and displayed isosurfaces in Fig. 10 show similar large-scale evolution for THINC/QQ and QMMS. At $t = 0.1$, the curved fronts and lower deformed layer occupy nearly the same locations. A small difference appears in the upper part of the deformed high-gradient front at $t = 0.3$ and becomes more pronounced by $t = 0.5$, while the overall structure remains closely aligned. The localized discrepancy occurs near the intersecting boundaries, where only interior-cell information is used to estimate the least-squares gradient and Hessian on a one-sided stencil. This stencil geometry may affect the fitted coefficients of the quadratic field $P$ and the resulting nonlinear face averages, which THINC/QQ and QMMS evaluate differently. Small differences in these face values can then affect the TBV values, BVD candidate selection, and numerical flux. As the density contact is transported along the symmetry boundary, these differences can accumulate and become more visible at later times.

### 3.4. Reconstruction-stage computational cost

The measured times are normalized by the numbers of time steps, Runge–Kutta stages, cells, and reconstructed variables to enable a more direct comparison of the computational cost per cell-variable reconstruction. All timings are obtained on a single physical CPU core. Sigmoid reconstruction is evaluated for every cell-variable pair during profiling. Table 2 summarizes the resulting reconstruction-stage timings together with the times required for surface-constant and face-average evaluations for each benchmark. These timings cover the reconstruction stage rather than the complete simulation.

Table 2. Reconstruction-stage computational timings and corresponding speedups of QMMS relative to THINC/QQ for the nine benchmark problems.

| (a) Reconstruction stage | | | |
|---|---|---|---|
| Case | THINC/QQ (ns) | QMMS (ns) | Speedup |
| Rotation | 1033 | 502 | 2.06× |
| Configuration 3 | 1036 | 477 | 2.17× |
| Shock–vortex | 1118 | 517 | 2.16× |
| Shock–mixing layer | 1113 | 514 | 2.16× |
| Mach 3 step | 1122 | 517 | 2.17× |
| Double Mach reflection | 1051 | 522 | 2.01× |
| 2D mean | 1079 | 508 | 2.12× |
| Deformation | 2418 | 1372 | 1.76× |
| Shock–cylinder | 2328 | 1205 | 1.93× |
| Spherical Riemann | 2284 | 1216 | 1.88× |
| 3D mean | 2343 | 1264 | 1.85× |

| (b) Surface-constant evaluation | | | |
|---|---|---|---|
| Case | THINC/QQ (ns) | QMMS (ns) | Speedup |
| Rotation | 321 | 37.6 | 8.5× |
| Configuration 3 | 347 | 36.9 | 9.4× |
| Shock–vortex | 380 | 36.5 | 10.4× |
| Shock–mixing layer | 390 | 36.8 | 10.6× |
| Mach 3 step | 387 | 35.3 | 10.9× |
| Double Mach reflection | 317 | 37.6 | 8.4× |
| 2D mean | 357 | 36.8 | 9.7× |
| Deformation | 380 | 25.1 | 15.2× |
| Shock–cylinder | 438 | 25.0 | 17.5× |
| Spherical Riemann | 430 | 23.9 | 18.0× |
| 3D mean | 416 | 24.6 | 16.9× |

| (c) Face-average evaluation | | | |
|---|---|---|---|
| Case | THINC/QQ (ns) | QMMS (ns) | Speedup |
| Rotation | 564 | 268 | 2.10× |
| Configuration 3 | 560 | 269 | 2.08× |
| Shock–vortex | 586 | 289 | 2.03× |
| Shock–mixing layer | 562 | 279 | 2.01× |
| Mach 3 step | 581 | 290 | 2.00× |
| Double Mach reflection | 554 | 267 | 2.07× |
| 2D mean | 568 | 277 | 2.05× |
| Deformation | 910 | 223 | 4.08× |
| Shock–cylinder | 940 | 220 | 4.26× |
| Spherical Riemann | 916 | 213 | 4.30× |
| 3D mean | 922 | 219 | 4.21× |

Across all nine benchmarks, QMMS has a lower reconstruction-stage computational cost than THINC/QQ. For the complete reconstruction stage in Table 2(a), the two-dimensional timings are 477–522 ns for QMMS and 1033–1122 ns for THINC/QQ, corresponding to case-wise speedups of 2.01–2.17. The corresponding three-dimensional times are 1205–1372 ns and 2284–2418 ns, with speedups of 1.76–1.93. The ratios of the mean THINC/QQ and QMMS times are 2.12 and 1.85 in two and three dimensions, respectively. For the components reported in Table 2(b) and (c), the mean-time ratios for surface-constant evaluation are 9.7 and 16.9 in two and three dimensions, respectively, while those for face-average evaluation are 2.05 and 4.21. Surface-constant evaluation has the larger speedup in both dimensions because its algebraic approximation replaces the quadrature sums within Newton iteration. The mean-row speedups are ratios of mean times rather than averages of case-wise speedups.

Although the component speedups are larger in three dimensions, the speedup of the complete reconstruction stage is lower. The two modified components account for approximately 85.7% of the baseline reconstruction time in two dimensions and 57.1% in three dimensions, based on the reported mean times. The remaining operations, including stencil-based polynomial reconstruction, limiting, and BVD selection, therefore occupy a larger fraction of the three-dimensional cost. This smaller accelerated fraction limits the overall reconstruction-stage speedup.

## 4. Conclusion

QMMS approximates the cell and face integrals in THINC/QQ reconstruction using the mean and variance of the quadratic field. A slope-matched Gaussian-CDF surrogate and a moment-matched Gaussian distribution provide closed-form approximations for the surface constant and face averages, eliminating Newton iteration and runtime quadrature in these evaluations. The hyperbolic-tangent profile, quadratic representation, candidate steepness parameters, and BVD selection procedure are otherwise unchanged.

Single-cell tests on six cell types show overlapping cell-average RMS error ranges for THINC/QQ and QMMS over the reported steepness interval. THINC/QQ generally gives lower median cell-average errors at small scaled-field magnitudes, whereas the QMMS error increases more gradually as the scaled-field magnitude increases. The face-average RMS ranges extend to larger upper values for QMMS. The flow benchmarks show that the principal transported and shock structures are largely preserved, while local differences appear mainly near discontinuities, contact regions, and rolled shear layers. In rigid-body rotation, the global $L_1$ error increases by approximately 2% with QMMS, although the object-wise errors vary with the transported profile. In the three-dimensional time-reversed deformation problem, the return-time $L_1$ error is approximately 4% higher for QMMS. The compressible-flow cases similarly show close agreement in the principal shock locations, with more visible local differences in secondary waves and shear-dominated regions.

The computational advantage is concentrated in the two integral-evaluation steps modified by QMMS. The mean THINC/QQ-to-QMMS time ratios for surface-constant evaluation are 9.7 in two dimensions and 16.9 in three dimensions, while those for face-average evaluation are 2.05 and 4.21. Under the profiling workload, the corresponding ratios of mean times for the complete reconstruction stage are 2.12 and 1.85. These results show that the moment-based approximation reduces the computational cost of THINC/QQ reconstruction and removes topology-dependent runtime quadrature from the modified integral evaluations, while the principal transported and shock structures remain closely aligned across the tested two- and three-dimensional benchmarks.

## Acknowledgements

The author has no acknowledgements to declare.

## CRediT authorship contribution statement

Young-Lin Yoo: Conceptualization, Methodology, Software, Validation, Formal analysis, Investigation, Visualization, Writing – original draft, Writing – review & editing.

## Declaration of competing interest

The author declares no known competing financial interests or personal relationships that could have influenced the work reported in this paper.

## Data availability

Data and code will be made available on reasonable request.

## Funding

This research did not receive any specific grant from funding agencies in the public, commercial, or not-for-profit sectors.

**Declaration of generative AI and AI-assisted technologies in the manuscript preparation process**

During the preparation of this work, the author used OpenAI ChatGPT to improve English grammar, sentence structure, and readability. After using this tool, the author reviewed and edited the content as needed and takes full responsibility for the content of the published article.

## References

[1] F. Xiao, Y. Honma, T. Kono, A simple algebraic interface capturing scheme using hyperbolic tangent function, International Journal for Numerical Methods in Fluids 48 (2005) 1023–1040. https://doi.org/10.1002/fld.975.

[2] K. Yokoi, Efficient implementation of THINC scheme: A simple and practical smoothed VOF algorithm, Journal of Computational Physics 226 (2007) 1985–2002. https://doi.org/10.1016/j.jcp.2007.06.020.

[3] F. Xiao, S. Ii, C. Chen, Revisit to the THINC scheme: A simple algebraic VOF algorithm, Journal of Computational Physics 230 (2011) 7086–7092. https://doi.org/10.1016/j.jcp.2011.06.012.

[4] S. Ii, K. Sugiyama, S. Takeuchi, S. Takagi, Y. Matsumoto, F. Xiao, An interface capturing method with a continuous function: The THINC method with multi-dimensional reconstruction, Journal of Computational Physics 231 (2012) 2328–2358. https://doi.org/10.1016/j.jcp.2011.11.038.

[5] Z. Sun, S. Inaba, F. Xiao, Boundary Variation Diminishing (BVD) reconstruction: A new approach to improve Godunov schemes, Journal of Computational Physics 322 (2016) 309–325. https://doi.org/10.1016/j.jcp.2016.06.051.

[6] X. Deng, Y. Shimizu, F. Xiao, A fifth-order shock capturing scheme with two-stage boundary variation diminishing algorithm, Journal of Computational Physics 386 (2019) 323–349. https://doi.org/10.1016/j.jcp.2019.02.024.

[7] X. Deng, Y. Shimizu, B. Xie, F. Xiao, Constructing higher order discontinuity-capturing schemes with upwind-biased interpolations and boundary variation diminishing algorithm, Computers & Fluids 200 (2020) 104433. https://doi.org/10.1016/j.compfluid.2020.104433.

[8] S. Takagi, L. Fu, H. Wakimura, F. Xiao, A novel high-order low-dissipation TENO-THINC scheme for hyperbolic conservation laws, Journal of Computational Physics 452 (2022) 110899. https://doi.org/10.1016/j.jcp.2021.110899.

[9] W. Zhang, N. Fleischmann, S. Adami, N.A. Adams, A hybrid WENO5IS-THINC reconstruction scheme for compressible multiphase flows, Journal of Computational Physics 498 (2024) 112672. https://doi.org/10.1016/j.jcp.2023.112672.

[10] F. Zeng, Y. Qiu, X. Lyu, Z. Jiang, W. Chen, High-fidelity WENO5IS-THINC-BVD hybrid reconstruction for compressible multiphase flows, Computers & Fluids 310 (2026) 107034. https://doi.org/10.1016/j.compfluid.2026.107034.

[11] K.-M. Shyue, F. Xiao, An Eulerian interface sharpening algorithm for compressible two-phase flow: The algebraic THINC approach, Journal of Computational Physics 268 (2014) 326–354. https://doi.org/10.1016/j.jcp.2014.03.010.

[12] X. Deng, S. Inaba, B. Xie, K.-M. Shyue, F. Xiao, High fidelity discontinuity-resolving reconstruction for compressible multiphase flows with moving interfaces, Journal of Computational Physics 371 (2018) 945–966. https://doi.org/10.1016/j.jcp.2018.03.036.

[13] S. Ii, B. Xie, F. Xiao, An interface capturing method with a continuous function: The THINC method on unstructured triangular and tetrahedral meshes, Journal of Computational Physics 259 (2014) 260–269. https://doi.org/10.1016/j.jcp.2013.11.034.

[14] B. Xie, S. Ii, F. Xiao, An efficient and accurate algebraic interface capturing method for unstructured grids in 2 and 3 dimensions: The THINC method with quadratic surface representation, International Journal for Numerical Methods in Fluids 76 (2014) 1025–1042. https://doi.org/10.1002/fld.3968.

[15] B. Xie, F. Xiao, Toward efficient and accurate interface capturing on arbitrary hybrid unstructured grids: The THINC method with quadratic surface representation and Gaussian quadrature, Journal of Computational Physics 349 (2017) 415–440. https://doi.org/10.1016/j.jcp.2017.08.028.

[16] D. Chen, B. Xie, F. Xiao, Revisit to the THINC/QQ scheme: Recent progress to improve accuracy and robustness, International Journal for Numerical Methods in Fluids 94 (2022) 719–755. https://doi.org/10.1002/fld.5072.

[17] L. Cheng, X. Deng, B. Xie, Y. Jiang, F. Xiao, Low-dissipation BVD schemes for single and multi-phase compressible flows on unstructured grids, Journal of Computational Physics 428 (2021) 110088. https://doi.org/10.1016/j.jcp.2020.110088.

[18] B. Xie, P. Jin, F. Xiao, An unstructured-grid numerical model for interfacial multiphase fluids based on multi-moment finite volume formulation and THINC method, International Journal of Multiphase Flow 89 (2017) 375–398. https://doi.org/10.1016/j.ijmultiphaseflow.2016.10.016.

[19] L. Cheng, X. Deng, B. Xie, Y. Jiang, F. Xiao, A new 3D OpenFoam solver with improved resolution for hyperbolic systems on hybrid unstructured grids, Applied Mathematical Modelling 108 (2022) 142–166. https://doi.org/10.1016/j.apm.2022.03.022.

[20] D. Chen, X. Tong, B. Xie, F. Xiao, Y. Li, An accurate and efficient multiphase solver based on THINC scheme and adaptive mesh refinement, International Journal of Multiphase Flow 162 (2023) 104409. https://doi.org/10.1016/j.ijmultiphaseflow.2023.104409.

[21] P. Jin, B. Xie, A multi-moment finite volume formulation for the interaction between free surface flow and moving bodies with THINC method, Computers & Fluids 265 (2023) 105994. https://doi.org/10.1016/j.compfluid.2023.105994.

[22] X. Tong, D. Chen, L. Cheng, B. Xie, A THINC-based numerical model for incompressible flows with free surfaces on the overset grids: On preserving accuracy and conservation of volume fraction, Journal of Fluids and Structures 135 (2025) 104294. https://doi.org/10.1016/j.jfluidstructs.2025.104294.

[23] H. Wakimura, T. Aoki, F. Xiao, A low-dissipation numerical method based on boundary variation diminishing principle for compressible gas–liquid two-phase flows with phase change on unstructured grid, Physics of Fluids 37 (2025) 016103. https://doi.org/10.1063/5.0243965.

[24] A.K. Pandare, J. Waltz, J. Bakosi, Multi-material hydrodynamics with algebraic sharp interface capturing, Computers & Fluids 215 (2021) 104804. https://doi.org/10.1016/j.compfluid.2020.104804.

[25] L. Qian, Y. Wei, F. Xiao, Coupled THINC and level set method: A conservative interface capturing scheme with high-order surface representations, Journal of Computational Physics 373 (2018) 284–303. https://doi.org/10.1016/j.jcp.2018.06.074.

[26] R. Kumar, L. Cheng, Y. Xiong, B. Xie, R. Abgrall, F. Xiao, THINC scaling method that bridges VOF and level set schemes, Journal of Computational Physics 436 (2021) 110323. https://doi.org/10.1016/j.jcp.2021.110323.

[27] Y. Xiong, F. Xiao, B. Xie, A hybrid volume of fluid and level set interface capturing scheme with quartic surface representation for unstructured meshes, International Journal for Numerical Methods in Fluids 94 (2022) 1542–1565. https://doi.org/10.1002/fld.5103.

[28] H.Y. Zhao, P.J. Ming, W.P. Zhang, J.K. Chen, A direct time-integral THINC scheme for sharp interfaces, Journal of Computational Physics 393 (2019) 139–161. https://doi.org/10.1016/j.jcp.2019.05.011.

[29] M. Huang, L. Cheng, W. Ying, X. Deng, F. Xiao, A low-dissipation reconstruction scheme for compressible single- and multi-phase flows based on artificial neural networks, Journal of Computational Physics 530 (2025) 113894. https://doi.org/10.1016/j.jcp.2025.113894.

[30] D. Chen, S. Liao, B. Xie, An accurate THINC scheme for interface capturing based on homotopy analysis method, Journal of Computational Physics 544 (2026) 114420. https://doi.org/10.1016/j.jcp.2025.114420.

[31] K. Zhang, Y. Shen, Interface capturing schemes based on sigmoid functions, Computers & Fluids 280 (2024) 106352. https://doi.org/10.1016/j.compfluid.2024.106352.

[32] D. Kim, C.B. Ivey, F.E. Ham, L.G. Bravo, An efficient high-resolution Volume-of-Fluid method with low numerical diffusion on unstructured grids, Journal of Computational Physics 446 (2021) 110606. https://doi.org/10.1016/j.jcp.2021.110606.

[33] W. Ni, Q. Zeng, Y. Ruan, Z. He, A novel steepness-adjustable harmonic volume-of-fluid method for interface capturing, Journal of Computational Physics 501 (2024) 112765. https://doi.org/10.1016/j.jcp.2024.112765.

[34] Z. He, Y. Ruan, Y. Yu, B. Tian, F. Xiao, Self-adjusting steepness-based schemes that preserve discontinuous structures in compressible flows, Journal of Computational Physics 463 (2022) 111268. https://doi.org/10.1016/j.jcp.2022.111268.

[35] T. Fukuda, S. Yamashita, H. Yoshida, Development of a new THINC/WLIC method based on a separate evaluation of the geometrical fidelity and the interface sharpness, Journal of Computational Physics 545 (2026) 114485. https://doi.org/10.1016/j.jcp.2025.114485.

[36] J. López, Unsplit geometric volume-of-fluid method with iterative piecewise-paraboloid interface reconstruction on arbitrary three-dimensional grids, Journal of Computational Physics 553 (2026) 114714. https://doi.org/10.1016/j.jcp.2026.114714.

[37] R.J. Spiteri, S.J. Ruuth, A New Class of Optimal High-Order Strong-Stability-Preserving Time Discretization Methods, SIAM J. Numer. Anal. 40 (2002) 469–491. https://doi.org/10.1137/S0036142901389025.

[38] T. Barth, D. Jespersen, The design and application of upwind schemes on unstructured meshes, in: 27th Aerospace Sciences Meeting, American Institute of Aeronautics and Astronautics, Reno, NV, USA, 1989. https://doi.org/10.2514/6.1989-366.

[39] K.H. Kim, C. Kim, Accurate, efficient and monotonic numerical methods for multi-dimensional compressible flows, Journal of Computational Physics 208 (2005) 570–615. https://doi.org/10.1016/j.jcp.2005.02.022.

[40] J.S. Park, S.-H. Yoon, C. Kim, Multi-dimensional limiting process for hyperbolic conservation laws on unstructured grids, Journal of Computational Physics 229 (2010) 788–812. https://doi.org/10.1016/j.jcp.2009.10.011.

[41] M.G. Duffy, Quadrature Over a Pyramid or Cube of Integrands with a Singularity at a Vertex, SIAM J. Numer. Anal. 19 (1982) 1260–1262. https://doi.org/10.1137/0719090.

[42] M. Gandhi, K. Lee, Y. Pan, E. Theodorou, Propagating Uncertainty through the tanh Function with Application to Reservoir Computing, arXiv preprint arXiv:1806.09431 (2018). https://doi.org/10.48550/arXiv.1806.09431.

[43] D.J.C. MacKay, The Evidence Framework Applied to Classification Networks, Neural Computation 4 (1992) 720–736. https://doi.org/10.1162/neco.1992.4.5.720.

[44] C.E. Rasmussen, C.K.I. Williams, Gaussian processes for machine learning, 3rd printing, MIT Press, Cambridge, Mass., 2008.

[45] R.J. LeVeque, High-Resolution Conservative Algorithms for Advection in Incompressible Flow, SIAM J. Numer. Anal. 33 (1996) 627–665. https://doi.org/10.1137/0733033.

[46] P.D. Lax, X.-D. Liu, Solution of Two-Dimensional Riemann Problems of Gas Dynamics by Positive Schemes, SIAM J. Sci. Comput. 19 (1998) 319–340. https://doi.org/10.1137/S1064827595291819.

[47] C.W. Schulz-Rinne, J.P. Collins, H.M. Glaz, Numerical Solution of the Riemann Problem for Two-Dimensional Gas Dynamics, SIAM J. Sci. Comput. 14 (1993) 1394–1414. https://doi.org/10.1137/0914082.

[48] G.-S. Jiang, C.-W. Shu, Efficient Implementation of Weighted ENO Schemes, Journal of Computational Physics 126 (1996) 202–228. https://doi.org/10.1006/jcph.1996.0130.

[49] K. Huang, H. Wu, H. Yu, D. Yan, Cures for numerical shock instability in HLLC solver, International Journal for Numerical Methods in Fluids 65 (2011) 1026–1038. https://doi.org/10.1002/fld.2217.

[50] H.C. Yee, N.D. Sandham, M.J. Djomehri, Low-Dissipative High-Order Shock-Capturing Methods Using Characteristic-Based Filters, Journal of Computational Physics 150 (1999) 199–238. https://doi.org/10.1006/jcph.1998.6177.

[51] P. Woodward, P. Colella, The numerical simulation of two-dimensional fluid flow with strong shocks, Journal of Computational Physics 54 (1984) 115–173. https://doi.org/10.1016/0021-9991(84)90142-6.

[52] A. Harten, P.D. Lax, B. van Leer, On Upstream Differencing and Godunov-Type Schemes for Hyperbolic Conservation Laws, SIAM Review 25 (1983) 35–61. https://doi.org/10.1137/1025002.

[53] D. Enright, R. Fedkiw, J. Ferziger, I. Mitchell, A Hybrid Particle Level Set Method for Improved Interface Capturing, Journal of Computational Physics 183 (2002) 83–116. https://doi.org/10.1006/jcph.2002.7166.

[54] J.O. Langseth, R.J. LeVeque, A Wave Propagation Method for Three-Dimensional Hyperbolic Conservation Laws, Journal of Computational Physics 165 (2000) 126–166. https://doi.org/10.1006/jcph.2000.6606.

[55] F. Alauzet, M. Mehrenberger, P1-conservative solution interpolation on unstructured triangular meshes, International Journal for Numerical Methods in Engineering 84 (2010) 1552–1588. https://doi.org/10.1002/nme.2951.

[56] T.F. Illenseer, W.J. Duschl, Two-dimensional central-upwind schemes for curvilinear grids and application to gas dynamics with angular momentum, Computer Physics Communications 180 (2009) 2283–2302. https://doi.org/10.1016/j.cpc.2009.07.016.

[57] G.-H. Tu, X.-J. Yuan, A characteristic-based shock-capturing scheme for hyperbolic problems, Journal of Computational Physics 225 (2007) 2083–2097. https://doi.org/10.1016/j.jcp.2007.03.007.

[58] E.F. Toro, M. Spruce, W. Speares, Restoration of the contact surface in the HLL-Riemann solver, Shock Waves 4 (1994) 25–34. https://doi.org/10.1007/BF01414629.

[59] X. Huang, J. Chen, J. Zhang, L. Wang, Y. Wang, An Adaptive Mesh Refinement–Rotated Lattice Boltzmann Flux Solver for Numerical Simulation of Two and Three-Dimensional

Compressible Flows with Complex Shock Structures, Symmetry 15 (2023) 1909. https://doi.org/10.3390/sym15101909.